\documentclass{article}
\usepackage{metric-min}
\hypersetup{pdftitle={Metric-minimizing surfaces revisited},
pdfauthor={Anton Petrunin and Stephan Stadler}}

\begin{document}

\title{Metric-minimizing surfaces revisited}
\author{Anton Petrunin and Stephan Stadler}
\nofootnote{A.~Petrunin was partially supported by NSF grant DMS 1309340.
S.~Stadler was supported by DFG grants STA 1511/1-1 and SPP 2026.}

\newcommand{\Addresses}{{\bigskip\footnotesize
Anton Petrunin, \par\nopagebreak\textsc{Department of Mathematics, PSU, University Park, PA 16802, USA}
\par\nopagebreak
\textit{Email}: \texttt{petrunin@math.psu.edu}

\medskip
 
Stephan Stadler,
\par\nopagebreak\textsc{Mathematisches Institut der Universit\"at M\"unchen, Theresienstr. 39, D-80333 M\"unchen, Germany}
\par\nopagebreak
\textit{Email}: \texttt{stadler@math.lmu.de}
}}

\date{}


\maketitle

\begin{abstract}
A surface that does not admit a length-nonincreasing deformation is called \emph{metric-minimizing}.
We show that metric-minimizing surfaces in $\CAT(0)$ spaces are locally $\CAT(0)$ with respect to their length metrics. 
\end{abstract}

\section{Introduction}

\parbf{Main result.}
Assume $s$ is a Lipschitz embedding of the disc $\DD$ into the Euclidean space $\RR^3$.
We say $s$ is \emph{metric-minimizing} if its intrinsic metric is minimal. Explicitly, this means that if a map $s'\:\DD\to \RR^3$ agrees with 
$s$ on $\partial \DD=\SS^1$ and fulfills 
\[\length s\circ\gamma\ge \length s'\circ\gamma\]
for all curves $\gamma$ in $\DD$, then equality holds for all curves $\gamma$.

As follows from the main theorem below, the induced metric on the disc for a metric-minimizing embedding is $\CAT(0)$.

To formulate the theorem in full generality we need to extend the definition of metric-minimizing maps in two ways.
First, we will only assume continuity of $s$, so the map $s$ might not be an embedding nor does it have to be Lipschitz;
this part is tricky --- a straightforward generalization turns out to be too weak.
Secondly, we will have to allow arbitrary metric spaces as targets; this part is straightforward. 

Let $Y$ be a metric space
and $s\:\DD\to Y$ be a continuous map.
Let $\gamma\:\SS^1\to Y$ be a closed rectifiable curve.
We say that $s\:\DD\to Y$ \emph{spans} $\gamma$,
if $s$ is an extension of $\gamma$; that is,  $s|_{\SS^1}=\gamma$.

Consider the induced length-pseudometric on $\DD$ defined as 
\[\<x-y\>_s=\inf_\alpha\{\,\length (s\circ\alpha)\,\},\]
where the greatest lower bound is taken over all paths $\alpha$ from $x$ to $y$ in $\DD$.
The distance $\<x-y\>_s$ might take infinite values.
We denote by $\<\DD\>_s$ the corresponding metric space;
see the next section for a precise definition.

The space $\<\DD\>_s$ comes with the projections $\DD\xrightarrow{\hat\pi_s}\<\DD\>_s\xrightarrow{\hat s}Y$;
the restriction $\hat\pi_s|_{\SS^1}$ will be denoted by $\delta_s$.

Assume that $s$ and $s'$ are two maps spanning the same curve.
We write $s\succcurlyeq s'$ if there is a \emph{majorization}
$\mu\:\<\DD\>_s\to \<\DD\>_{s'}$;
meaning that $\mu$ is a \emph{short} map and $\delta_{s'}=\mu\circ\delta_s$;
here and further \emph{short} stands for \emph{distance-nonincreasing}.

(Note that if $s$ and $s'$ are Lipschitz embeddings, then any majorization $\<\DD\>_s\to \<\DD\>_{s'}$ admits a continuous lifting $\DD\to\DD$; 
in general, such a lifting may not exist.)

A map $s\:\DD\to Y$ will be called a \emph{metric-minimizing disc} if $\delta_s$ is rectifiable and $s\succcurlyeq s'$ implies that the corresponding majorization $\mu$ is distance-preserving.\label{metric-minimizing disc}

A topological space $W$ together with a choice of a closed curve $\delta\:\SS^1\zz\to W$ is called a \emph{disc-retract} 
if the mapping cylinder of~$\delta$ 
\[W_\delta=W\bigsqcup_{\delta(u)\sim(u,0)}\SS^1\times[0,1]\]
is homeomorphic to the disc $\DD$.
The curve $\delta$ will be called the \emph{boundary curve} of the disc-retract $W$.

\begin{wrapfigure}{r}{28 mm}
\begin{lpic}[t(-2 mm),b(-0 mm),r(0 mm),l(0 mm)]{pics/two-bry-curves(1)}
\end{lpic}
\end{wrapfigure}

An example of a disc-retract that is not a disc is shown in the picture;
note that it has essentially different boundary curves,
in the sense that one is not a reparametrization of the other.
The definition is motivated by the following observation: \textit{if the boundary curve of a disc-retract $W$ is a simple closed curve, then $W$ is homeomorphic to $\DD$.}

\begin{thm}{Main theorem}\label{thm:main}
Assume $Y$ is a $\CAT(0)$ space and $s\:\DD\to Y$ is a metric-minimizing disc.
Then $\<\DD\>_s$ is a $\CAT(0)$ disc-retract with boundary curve $\delta_s$.

In particular, if $s$ spans a simple closed curve, then $\<\DD\>_s$ is a $\CAT(0)$ disc.
\end{thm}

Note that apart from continuity we did not make any regularity assumptions on the map $s$;
in particular, beforehand, the space $\<\DD\>_s$ might have wild topology.

If we remove the condition that the boundary curve is rectifiable, then the space $\<\DD\>_s$ might have points at infinite distance from each other.
However, our proof shows that all triangles with finite sides in $\<\DD\>_s$ are still thin; 
in particular, each metric component of $\<\DD\>_s$ is a $\CAT(0)$ space.
So, one could consider $\<\DD\>_s$ as a $\CAT(0)$ space where infinite distances between points are legal.

The class of metric-minimizing discs is huge --- if the ambient space $Y$ is $\CAT(0)$, then one can find a metric-minimizing  disc $\succcurlyeq$-\emph{below} for
any map $s\:\DD\zz\to Y$ such that the metric $(x,y)\mapsto\<x,y\>_s$ is continuous
(it is an application of the ultralimit+projection construction described below).
Metric-minimizing discs include many well-studied maps from $\DD$ to metric spaces.
The ruled discs considered by Alexandrov \cite{A} are evidently metric-minimizing since one cannot shorten a geodesic. 
In addition, harmonic maps in the sense of Korevaar--Schoen \cite{KS} from the disc to $\CAT(0)$ spaces are metric-minimizing.
Indeed, harmonic discs are solutions to the Dirichlet problem;
that is energy minimizing fillings of given loops.
Since the Korevaar--Schoen energy is convex, such solutions are unique. But decreasing the intrinsic metric will only decrease the energy of the corresponding map.
As a consequence, minimal discs in the sense of Lytchak--Wenger \cite{LW} are metric-minimizing;
these are Douglas--Rado solutions to the Plateau problem and therefore harmonic.

As intended, the main theorem subsumes and generalizes several previously known results, including Alexandrov's theorem about ruled surfaces in \cite{A} and 
the main theorem by the first author announced in \cite{petrunin-metric-min}, assuming it is formulated correctly; see Section~\ref{sec:old-thm} and \cite{petrunin-metric-min-correction}. 
It is closely related to the classical Gauss formula and results on saddle surfaces by Samuel Shefel \cite{shefel-2D,shefel-3D}.
It also generalizes the result on minimal surfaces by Alexander Lytchak and Stefan Wenger 
\cite[Theorem 1.2]{LW5}  and an earlier result of Chikako Mese \cite{mese};
see also \cite[Chapter 4]{akp} and \cite{petrunin-stadler}. 
Despite that some special cases of Theorem \ref{thm:main} were known, the result is new even for harmonic discs.

Let us list a few applications of metric-minimizing surfaces.
They were used by Alexander Lytchak \cite{L} to study sets of positive reach (ruled surfaces).
In  \cite{AB}, Stephanie Alexander and Richard Bishop used them to generalize the Gauss equation to non-smooth spaces (ruled surfaces).
In \cite{G}, Mikhael Gromov used them to bound the complexity of smooth maps (general metric-minimizing surfaces).
In \cite{LS}, Alexander Lytchak and the second author used them to deform general $\CAT(0)$ spaces (minimal discs).
In \cite{St}, the second author used them in the proof of a $\CAT(0)$ version of the Fary--Milnor theorem to control the mapping behavior of minimal surfaces (minimal discs).  

Note that any smooth metric-minimizing disc in a Euclidean space is saddle.
The converse fails in general (see Section~\ref{sec:smooth}), but the following statement gives a local converse in dimension three.
In the four-dimensional case, we expect that even the local converse does not hold.

\begin{thm}{Proposition}\label{prop:smooth}
Any smooth strictly saddle surface in $\RR^3$ is locally metric-minimizing.
\end{thm}

\parbf{Structure of the paper.}
In Section~\ref{sec:defs} we define a couple of metrics induced by a given continuous map.
This naturally leads to the monotone-light factorization theorem in a metric context.
 
In Section~\ref{Metric-minimizing discs} we obtain topological control on the length space
associated to a metric-minimizing disc. 
This part follows from Moore's quotient theorem on cell-like maps.
 
Section~\ref{Compactness} establishes compactness of a certain class of nonpositively curved surfaces.
This will be used later, when we approximate metric-minimizing maps.
 
In Section~\ref{Key Lemma} we prove the key lemma which says that the restriction to a finite set of
a metric-minimizing disc factorizes as a composition of short maps over a nonpositively curved surface.
The proof uses the properties of \emph{metric-minimizing graphs} discussed in Section~\ref{Metric-minimizing graphs}.
 
In Section~\ref{Finite-whole extension lemma} we prove an extension result which is of independent interest and logically detached
from the rest of the paper.
It gives a criterion that allows to extend maps from subsets
to the whole space, once extensions on finite subsets are guaranteed.
 
In Section~\ref{Main theorem} we assemble the proof of the main theorem.
We show that a metric-minimizing disc factorizes as a composition of short maps over a
nonpositively curved surface.
The proof is finished by showing that the first map of this factorization
induces an isometry.

In Section~\ref{sec:old-thm}, we provide rectified formulations of the main result in \cite{petrunin-metric-min} and 
show that they follow from the main theorem.

In the last section, we discuss the relation between saddle surfaces and metric-minimizing discs.

\parbf{Acknowledgement.}
We want to thank 
Sergei Ivanov, 
Nina Lebedeva,
Carlo Sinestrari, 
Peter Topping 
and Burkhard Wilking 
for help.
We also would like to thank Alexander Lytchak for explaining his recent work with Stefan Wenger and for several helpful discussions.

\section{Definitions}
\label{sec:defs}

\parbf{Metrics and pseudometrics.}
Let $X$ be a set.
A \emph{pseudometric} on $X$ 
is a function $X\times X\to[0,\infty]$ denoted as $(x,y)\mapsto |x-y|$
such that 
\begin{itemize}
\item $|x-x|=0$, for any $x\in X$;
\item $|x-y|=|y-x|$, for any $x,y\in X$;
\item $|x-y|+|y-z|\ge|x-z|$ for any  $x,y,z\in X$.
\end{itemize}
If in addition $|x-y|=0$ implies $x=y$, then the pseudometric $|{*}-{*}|$ is called a \emph{metric}; some authors prefer to call it an \emph{$\infty$-metric} to emphasize that the distance between points might be infinite.
The value $|x-y|$ will also be called the \emph{distance} form $x$ to $y$.

A (pseudo)metric space $X$ is the underlying set (also denoted by $X$) equipped with a (pseudo)metric which often will be denoted by $|{*}-{*}|_X$.
We will use $X$ as an index if we want to emphasize that we are working in the space $X$;
for example, the ball of radius $R$ centered at $z$ in $X$ can be denoted as
\[B(z,R)_X=\set{x\in X}{|z-x|_X<R}.\]

For any pseudometric on a set there is an equivalence relation ``$\sim$'' such that 
\[x\sim y\quad\iff\quad|x-y|=0.\]
The pseudometric induced  on the set of equivalence classes 
\[[x]=\set{x'\in X}{x'\sim x}\] becomes a metric.
The obtained metric space will be denoted as $[X]$;
it comes with the projection map $X\to [X]$ defined as $x\mapsto [x]$.

For a metric space we can consider the equivalence relation ``$\approx$'' defined as 
\[x\approx y\quad\iff\quad|x-y|<\infty.\]
Its equivalence classes are called \emph{metric components}.
Note that by definition each metric component is a \emph{genuine metric space}, meaning that distances between points are finite.
Consequently, any metric space is a disjoint union of genuine metric spaces.

\parbf{Pseudometrics induced by a map.}
Assume $X$ is a topological space and $Y$ is a metric space.
Let $f\:X\to Y$ be a continuous map. 

Let us define the \emph{length pseudometric} on $X$ induced by $f$ as
\[\< x-y\>_f
=
\inf
\set{\length(f\circ\gamma)_Y}{\gamma\ \text{a path in}\  X\ \text{from}\ x\ \text{to}\ y}.\]
Denote by $\< X\>_f$ the corresponding metric space; that is,
\[\< X\>_f=[(X,\<{*}-{*}\>_f)].\] 

Similarly, define a \emph{connecting pseudometric} $|{*}-{*}|_f$ on $X$ by
\[|x-y|_f=\inf\{\diam f(K)\},\]
where the greatest lower bound is taken over all connected sets $K\subset X$ that contain $x$ and $y$;
if there is no such set, then $|x-y|_f=\infty$. 
The associated metric space will be
denoted as $|X|_f$;
that is,
\[| X|_f=[(X,|{*}-{*}|_f)].\]

For the projections 
\[\bar \pi_f\:X\to |X|_f\quad\text{and}\quad \hat \pi_f\:X\to \<X\>_f\]
we will also use the shortcut notations 
\[\bar x=\bar\pi_f(x) \quad\text{and}\quad  \hat x= \hat \pi_f(x).\]

\begin{thm}{Lemma}\label{lem:picont}
Let $X$ be a locally connected topological space and $Y$ a metric space. Assume that $f\:X\to Y$ is a continuous map. 
Then $\bar\pi_f\:X\to|X|_f$ is continuous.

In particular, if in addition $X$ is compact, then so is $|X|_f$.
\end{thm}

\parit{Proof.}
For a point $x\in X$ and $\eps>0$ we denote by $U$ the connected component of $f^{-1}[B(f(x),\eps)_{Y}]$ that contains $x$.
Since $X$ is locally connected, the set $U$ is open.

Note that $\bar\pi_f(U)\subset B(\bar x,2\cdot \eps)_{|X|_f}$;
hence the result.
\qeds

\begin{wrapfigure}{R}{37 mm}
\begin{tikzpicture}[scale=1.5]

  \node (1) at (1,.5) {$|X|_f$};
  \node (2) at (1,1.5){$\<X\>_f$};
  \node (11) at (2,2){$Y$};
  \node (12) at (0,2) {$X$};
\draw[
    >=latex,
    auto=right,                      
    loop above/.style={out=75,in=105,loop},
    every loop,
    ]
   (2) edge node{$\tau_f$}(1)
   (12) edge[bend left] node[swap]{$f$}(11)
   (12) edge node[swap]{$\hat\pi_f$}(2)
   (2) edge node[swap]{$ \hat f$}(11)
   (12) edge[bend right] node{$\bar \pi_f$}(1)
   (1) edge[bend right] node{$\bar f$}(11);
\end{tikzpicture}
\label{diagram-page}
\end{wrapfigure}

Note that $\tau_f\: \hat x\to \bar x$ defines a map $\tau_f\: \<X\>_f\zz\to |X|_f$, and by construction, it preserves the lengths of all curves coming from $X$.
Since $\<X\>_f$ is a length space, it implies that $\tau_f$ is short.
The map $\tau_f$ might not induce an isometry
\[\<X\>_f\to\<|X|_f\>_{\bar f}.\]
Moreover, $\tau_f$ does not have to be injective, an example is given in \cite[4.2]{petrunin-intrinisic}.
However, for metric-minimizing discs $f\:\DD\to Y$ both statements hold true; see
Proposition~\ref{prop:|D|}.

The space $\<|X|_f\>_{\bar f}$ is the intrinsic metric on $|X|_f$.
We will denote it briefly by $\<|X|\>_f$.
The corresponding pseudometric will be called \emph{intrinsic pseudometric} on $X$ induced by $f$; it will be denoted by $\<|{*}-{*}|\>_f$.
This is a more natural way to pullback intrinsic metric to $X$.
If $X$ is compact, then $\<|{*}-{*}|\>_f$ coincides with the pseudometric $\mathrm{pull}_f$ defined in \cite{petrunin-intrinisic}.
It will show up in Section~\ref{sec:old-thm}.

The maps $\bar f\:|X|_f\to Y$ and $ \hat f\:\<X\>_f\to Y$ are uniquely defined by the identity
\[f(x)=\bar f(\bar x)=  \hat f( \hat x)\] for any $x\in X$.
By construction, the diagram commutes.

Moreover, if $X$ is compact, then $\bar\pi_f$ has connected fibers, see Lemma \ref{cor:fiberconnected}.

\begin{thm}{Lemma}\label{lem:lengthpres}
Let $X$ be a compact metric space. 
Then $\bar f\:|X|_f\to Y$ preserves the length of every curve.
\end{thm}

\parit{Proof.}
Since $\bar f$ is short, we have $\length(\bar f\circ\gamma)\leq \length(\gamma)$ for every curve $\gamma$ in $|X|_f$.
Let a rectifiable curve $\gamma$ be given and let $\eps>0$. Choose points $x_i$ on $\gamma$ such that $\length(\gamma)\leq\sum_i |x_i,x_{i+1}|_{|X|_f}+\eps$.
Denote by $\gamma_i$ the piece of $\gamma$ between $x_i$ and $x_{i+1}$. Since $\bar\pi_f$ has connected fibers, each set $\bar\pi_f^{-1}(\gamma_i)$
is connected. Hence 
\begin{align*}
\length(\gamma)&\leq \sum_i \diam f(\bar\pi_f^{-1}(\gamma_i))+\eps=
\\
&=\sum_i \diam(\bar f(\gamma_i))+\eps\leq
\\
&\leq\length(\bar f\circ\gamma)+\eps.
\end{align*}
The claim follows since $\eps$ was arbitrary.
\qeds

\parbf{Metrics induced by metrics.}
If $X$ is a metric space, the two constructions above can be applied to the identity map $\id\:X\to X$.
In this case, the obtained spaces $\<X\>_\id$ and $|X|_\id$ will be denoted by $\<X\>$ and $|X|$ respectively.
The space $\<X\>$ is $X$ equipped with induced length metric.
All three spaces $\<X\>$, $|X|$ and $X$ have the same underlying set;
in other words, they can be considered as a single space with different metrics and tautological maps between them.
Both tautological maps 
\[\<X\>\to |X|\to X\]
are short and length-preserving.

As a consequence of Lemma~\ref{lem:lengthpres}, the tautological map $\<|X|_f\>_{\bar f}\zz\to \<|X|_f\>$ is an isometry;
that is, for any continuous map $f\:X\to Y$ the induced length metric on $|X|_f$ coincides with the length metric induced by 
$\bar f\:|X|_f\to Y$.

Recall that a geodesic in a metric space is a curve whose length coincides with the distance between its endpoints.
A metric space is called \emph{geodesic} if any two points at a finite distance can be joined by a geodesic.

\begin{thm}{Lemma}\label{lem:geospace}
Let $X$ be a compact metric space. 
Then $\<X\>$ is a complete geodesic space.
\end{thm}

The second statement is classical (see for example \cite[II-\S8 Thm. 3]{KF});
the first one appears in \cite{HK}.

\parit{Proof.}
Assume that $\<X\>$ is not complete.
Fix a Cauchy sequence $(x_n)$ in $\<X\>$ that is not converging in $\<X\>$.
After passing to a subsequence, we can assume that the points of the sequence appear on a rectifiable curve $\hat\gamma\:[0,1)\zz\to\<X\>$ in the same order.

The corresponding curve $\gamma\:[0,1)\to X$ has the same length.
Since $X$ is compact we can extend it to a path $\gamma_+\:[0,1]\to X$.
The curve 
\[\hat\gamma_+=\hat\pi\circ\gamma_+\:[0,1]\to\<X\>\]
has the same length.
Therefore $\hat\gamma_+(1)$ is the limit of $(x_n)$, a contradiction.

It remains to show that $\<X\>$ is geodesic.
Assume $\gamma_n$ is a sequence of constant speed paths from $x$ to $y$ in $X$
such that $\length(\hat\gamma_n)\to \<x-y\>$ as $n\to\infty$.
Since $X$ is compact, we can pass to a partial limit $\gamma$ of  $\gamma_n$.
The corresponding curve $\hat\gamma=\hat \pi\circ\gamma$ is the needed geodesic from $\hat x$ to $\hat y$ in $\<X\>$.
\qeds

\parbf{Monotone-light factorization.} 
Let $f\:X\to Y$ be a map between topological spaces.
Recall that 
\begin{itemize}
\item $f$ is called \emph{monotone} if the inverse image of each point is connected,
 \item $f$ is called \emph{light} if the inverse image of any point is totally disconnected.
\end{itemize}
Since a connected set has to be nonempty, \textit{any monotone map is onto}.

\begin{thm}{Lemma}\label{cor:fiberconnected}
Assume $X$ is a locally connected compact metric space and $Y$ is a metric space.
Let $f\:X\to Y$ be a continuous map.
Then the map $\bar \pi_f$ is monotone and $\bar f$ is light.
In particular, 
\[f=\bar f\circ\bar\pi_f\]
is a monotone-light factorization. 
\end{thm}

\parit{Proof.}
First, we prove the monotonicity of $\bar\pi_f$.

Assume the contrary;
that is, for some $x\in X$ the equivalence class 
\[K=\bar\pi_f^{-1}(\bar x)=\set{x'\in X}{|x-x'|_f=0}\]
is not connected. Since $X$ is normal, we
can cover $K$ by disjoint open sets $U,V\subset X$ such that both intersections
$K\cap U$ and $K\cap V$ are nonempty.

By Lemma~\ref{lem:picont}, $K$ is closed.

Suppose that $x\in U$ and pick $x'\in K\cap V$.
Then there is a sequence of connected sets $K_n\ni x,x'$ such that $\diam f(K_n)<\tfrac1n$.
For each $n$ we choose a point $k_n\in K_n\backslash (U\cup V)$.
Let $k$ be a partial limit of the sequence $(k_n)$.
It follows that $k\in K\backslash (U\cup V)$, a contradiction. 

Assume $\bar f$ is not light;
that is, the inverse image of some $y\in Y$ contains a closed connected set $C\subset |X|_f$ with more than one point.  
Note that the inverse image $Z:=\bar\pi_f^{-1}C$ is connected due to the monotonicity of $\bar\pi_f$. 
It follows that $|z-z'|_f=0$ for any two points $z,z'\in Z$, a contradiction.
\qeds

\section{Disc retracts}\label{Metric-minimizing discs}

A disc-retract as defined above is nothing but the image of strong deformation retraction of $\DD$;
the restriction of the retraction to the boundary can be taken as the corresponding boundary curve.
We will not need this statement, but it follows from Moore's quotient theorem quoted below. 

\begin{thm}{Proposition}\label{prop:|D|}
Let $Y$ be a metric space and $s\:\DD\to\ Y$ be a metric-minimizing disc.
Then $|\DD|_s$ is a disc-retract with boundary curve $\delta_s=\bar\pi_s|_{\SS^1}$.
Moreover, the map $\tau_s\:\<\DD\>_s\to |\DD|_s$ is injective and defines an isometry
$\<\DD\>_s\zz\to \<|\DD|\>_s$;
that is, $\tau_s$ is an isometry from $\<\DD\>_s$ to $|\DD|_s$ equipped with induced length metric.
\end{thm}

We need a little preparation before giving the proof.

Let $Y$ be a metric space and
$s\:\DD\to Y$ be a continuous map.
We say that $s$ has \label{page:no-bubble}\emph{no bubbles}
if for any point $p\in Y$ every connected component of the complement $\DD\backslash s^{-1}\{p\}$ contains a point from $\partial \DD$.

\begin{thm}{Lemma}\label{prop:point-complement}
Let $Y$ be a metric space and $s\:\DD\to Y$ be a metric-minimizing disc.
Then $s$ has no bubbles.
\end{thm}

\parit{Proof.}
Assume the contrary;
that is, there is $y\in Y$ such that the complement $\DD\backslash s^{-1}(y)$ contains a connected component $\Omega$ with $\partial \DD\cap \Omega=\emptyset$.

Let us define a new map $s'\:\DD\to\ Y$ by setting $s'(z)=y$ for any $x\in \Omega$ and $s'(x)=s(x)$ for any $x\notin \Omega$.

By construction, $s'$ and $s$ agree on $\partial\DD$. Moreover, $s\succcurlyeq s'$
because of the majorization $\mu\:\hat\pi_s(x)\zz\mapsto \hat\pi_{s'}(x)$.

Note that
\[\<x-x'\>_{s}>0=\<x-x'\>_{s'}\]
for a pair of distinct points $x,x'\in \Omega$.
In particular, $\mu$ is not distance-preserving, a contradiction.
\qeds

\begin{thm}{Lemma}\label{prop:disc-moore}
Let $Y$ be a metric space. 
Assume that a map $f\:\DD\to Y$ has no bubbles.
Then $|\DD|_f$ is homeomorphic to a disc-retract with boundary curve~$\bar\pi_f|_{\SS^1}$.
\end{thm}

This lemma is nearly identical to \cite[Corollary 7.12]{LW3}.
It could be considered a disc version of Moore's quotient theorem \cite{moore, daverman}.
The latter states that \textit{if a continuous map $f$ from the sphere $\SS^2$ to a Hausdorff space $X$
has acyclic fibers, then $f$ can be approximated by a homeomorphism};
in particular, $X$ is homeomorphic to $\SS^2$.

\parit{Proof.}
From Lemma~\ref{lem:picont} we know that $\bar\pi_f$ is continuous and hence $|\DD|_f$
is a compact metric space.

The mapping cone over $\DD$ along its boundary is homeomorphic to the sphere $\SS^2$;
denote by $\Sigma$ the mapping cone over $|\DD|_f$ with respect to $\bar\pi_s|_{\SS^1}$.
Let us extend the map $\bar\pi_f$ to a map between the mapping cones $\SS^2\to\Sigma$.
Note that this map satisfies Moore's quotient theorem, hence the statement follows.
\qeds

\parit{Proof of Proposition~\ref{prop:|D|}.}
The first two statements follow from Lemma~\ref{prop:point-complement} and Lemma~\ref{prop:disc-moore}.

Since $|\DD|_s$ is a disc-retract, the mapping cylinder over the boundary curve of $|\DD|_s$ is homeomorphic to $\DD$.
Denote by $r \:\DD\to |\DD|_s$ the corresponding retraction.

Note that $\<\DD\>_{\bar s\circ r}$ is isometric to $|\DD|_s$ equipped with the induced length metric.
Recall that the map $\tau_s\:\<\DD\>_s\to|\DD|_s$ is short and the induced map $\mu\:\<\DD\>_s\zz\to\<\DD\>_{\bar s\circ r}$ is a majorization.
Since $s$ is metric-minimizing, $\mu$ is distance-preserving.
Hence the statement follows.
\qeds

Assume $W$ is a disc-retract with a boundary curve $\delta$.
Recall that a point $p$ in a connected space $W$ is a \emph{cut point} if the complement $W\backslash\{p\}$ is not connected.

\begin{thm}{Lemma}\label{lem:discs}
Suppose that $W$ is a disc-retract.
Let $\Delta\subset W$ be a maximal connected subset that contains no cut points.
Assume $\Delta$ has at least two points.
Then the closure of $\Delta$ is homeomorphic to $\DD$.
\end{thm}

\parit{Proof.}
Let $\delta$ be a boundary curve of $W$, so the mapping cylinder $W_\delta$ is homeomorphic to $\DD$.

Note that $p$ is a cut point of $W$ if and only if $\delta^{-1}\{p\}$ has at least two connected components.
In particular, any cut point of $W$ lies on its boundary curve.

Denote by $\bar\Delta$ the closure of $\Delta$.
Note that for any $x\notin\bar\Delta$ there is a cut point $p\in\bar\Delta$ that cuts $x$ from $\Delta$.
Moreover, the map $\sigma\:x\mapsto p$ is uniquely defined on $W\backslash\bar\Delta$;
extend this map to the whole $W$ by identity in $\bar\Delta$.
By Moore's theorem, $\bar\Delta$ is a disc-retract with a boundary curve $\sigma\circ\delta$.
Namely, we apply Moore's theorem to $\SS^2= W_\delta/(\SS^1\times 1)$ and the quotient map $\SS^2\to \SS^2/\sim$ for the minimal equivalence relation such that $x\sim y$ if $\sigma(x)=\sigma(y)$.

The space $\bar\Delta$ has no cut points; in other words, $\sigma\circ\delta$ is monotonic.
It follows that $\sigma\circ\delta$ can be reparameterized into a simple closed curve.
By Jordan--Schoenflies theorem, the statement follows.
\qeds

The following lemma will be used in the final step in the proof of the main theorem, Section~\ref{Main theorem}.

\begin{thm}{Lemma}\label{lem:maj is isom}
Let $Y$ be a metric space and $s\:\DD\to Y$ be a metric-minimizing map.
Assume that there is a $\CAT(0)$ disc-retract $W$ with boundary curve $\delta$ and a short map $f\:\<\DD\>_s \to W$
such that $f\circ \delta_s=\delta$. If there exists a short map 
$q\: W\to Y$ with $q\circ \delta=s|_{\partial \DD}$, then the map $f$ is an isometry.
\end{thm}

\parit{Proof.}
Let $r\:\DD\to W$ be the projection from the mapping cylinder $\DD=W_\delta$. 
Note that $r$ is a retraction, $r|_{\partial \DD}=\delta$ and the composition $\DD\xrightarrow{r}W\xrightarrow{q} Y$ fulfills \[q\circ r|_{\partial \DD}=s|_{\partial \DD}.\]

Note that  $\<W\>_q=\<\DD\>_{q\circ r}$ and the natural projection $\rho\: W\to \<W\>_q$ is short.
It follows that $\rho\circ f\: \<\DD\>_s\to \<\DD\>_{q\circ r}$ is a majorization.
Since $s$ is metric-minimizing, $\rho\circ f$ is distance-preserving. 

Therefore $f$ is an isometric embedding. 
By Lemma \ref{lem:geospace} and Proposition \ref{prop:|D|}, $\<\DD\>_s$ is a complete geodesic space.
Therefore, the image of $f$ is a closed convex set in $W$, and it contains $\delta$.
Hence $f$ has to be surjective, and the lemma follows.
\qeds

\section{Compactness lemma}\label{Compactness}

A sequence of pairs $(X_n,\gamma_n)$, where $X_n$ is a metric space and $\gamma_n\:\SS^1\to X_n$ is a 
closed curve is said to \emph{converge} to $(X_\infty,\gamma_\infty)$ if there is a convergence of $X_n$ to $X_\infty$ 
in the sense of Gromov--Hausdorff for which $\gamma_n$ converges to $\gamma_\infty$ pointwise.

More precisely, we ask the following.
\begin{enumerate}[(1)]
\item There is a metric $\rho$ on the disjoint union 
\[\bm{X}=X_\infty\sqcup X_1\sqcup X_2\dots\]
that restricts to the given metric on each $X_\alpha$, for $\alpha\in\{1,2,\dots,\infty\}$, 
and such that $X_n$ converge to $X_\infty$ in the sense of Hausdorff as subsets in $(\bm{X},\rho)$.
\item  The sequence of compositions $\gamma_n\:\SS^1\to X_n \hookrightarrow\bm{X}$ 
converges to $\gamma_\infty\:\SS^1\zz\to X_\infty \hookrightarrow\bm{X}$ pointwise.
\end{enumerate}
Consider the class $\mathcal{K}_\ell$ of $\CAT(0)$ disc-retracts whose marked
boundary curves have Lipschitz constant $\ell$.

\begin{thm}{Compactness lemma}\label{lem:compact}
$\mathcal{K}_\ell$ is compact in the topology described above.
\end{thm}

The lemma follows from the two lemmas below.

\begin{thm}{Lemma}\label{lem:precompact}
$\mathcal{K}_\ell$ is precompact in the topology described above.
\end{thm}

\parit{Proof.}
Let $K$ be a metric space with the isometry class in $\mathcal {K}_\ell$.

Denote by $\area A$ the two-dimensional Hausdorff measure of $A\subset K$.
By the Euclidean isoperimetric inequality, we have 
\[\area K \le \pi\cdot\ell^2.\]

Fix $\eps>0$. 
Set $m=\lceil 10\cdot\tfrac\ell\eps\rceil$.
Choose $m$ points $y_1,\dots,y_m$ on $\partial K$ that divide $\partial K$ into arcs of equal length.

Consider the maximal set of points $\{x_1,\dots,x_n\}$ such that $d(x_i,x_j)>\eps$ and $d(x_i,y_j)>\eps$.

Note that the set $\{x_1,\dots,x_n,y_1,\dots,y_m\}$
is an $\eps$-net in $(K,d)$.
Further, note that the balls $B_i=B_{\eps/2}(x_i)$
do not overlap.

By comparison,
\[\area B_i\ge \tfrac{\pi\cdot\eps^2}{4}.\]
It follows that $n\le 4\cdot\left(\tfrac\ell\eps\right)^2$.
In particular, there is an integer-valued function $N(\eps)$, such that any  
$K$ as above contains an $\eps$-net
with at most $N(\eps)$ points.

The latter means that the class $\mathcal{K}_\ell$ is uniformly totally bounded.
By the selection theorem \cite[7.4.15]{BBI}, the class of metrics with this property is precompact in the Gromov--Hausdorff topology.

Since the set of $\ell$-Lipschitz maps defined on $\SS^1$ with compact target is compact 
with respect to pointwise convergence, we conclude that $\mathcal{K}_\ell$ is precompact in the topology defined above. 
\qeds

\begin{thm}{Lemma}\label{lem:closed}
$\mathcal{K}_\ell$ is closed in the topology described above.
\end{thm}

\parit{Proof.}
Let $(X_n,\gamma_n)$ be a sequence in $\mathcal{K}_\ell$.
Assume $X_n\to X$ and $\gamma_n\to\gamma$. 
Choose a point $o_n\in \partial X_n$ and define
$f_n:\DD\to X_n$ by sending the geodesic $[0,\theta]$ for $\theta\in\partial \DD=\SS^1$ to the geodesic path $[o_n,\gamma_n(\theta)]$ with constant speed. 

By comparison, $f_n$ is a $(5\cdot\ell)$-Lipschitz continuous ruled disc. 
The limit map $f\:\DD\to X$ is also a $(5\cdot\ell)$-Lipschitz continuous ruled disc.
In particular, $f$ is metric-minimizing.
By Lemma~\ref{prop:point-complement}, $f$ has no bubbles.
Therefore, by Lemma~\ref{prop:disc-moore}, $X$ is a disc-retract with boundary curve $\gamma$.
\qeds

\section{Metric-minimizing graphs}\label{Metric-minimizing graphs}

\emph{Metric-minimizing graphs} are defined analogously to metric-minimizing discs.

Namely, let $Y$ be a metric space, $\Gamma$ be a finite graph and $A$ be a subset of its vertexes.
Given two maps $f,f'\:\Gamma\to Y$, we write $f\succcurlyeq f'\rel A$ if $f$ and $f'$ agree on $A$ 
and there is a majorization $\mu\:\<\Gamma\>_f\to \<\Gamma\>_{f'}$
such that $f(a)\zz=f'(\mu(a))$ for any $a\in A$.

A map $f\:\Gamma\to Y$ is called \emph{metric-minimizing relative to $A$} if $f\succcurlyeq f'\rel A$ implies that the majorization $\mu$ is distance-preserving.

\begin{thm}{Proposition}\label{prop:metric-min-graph-exist}
Let $Y$ be a $\CAT(0)$ space, 
$\Gamma$ be a finite graph and $A$ be a subset of its vertexes.

Given a continuous map $f\:\Gamma\to Y$ there is a map $h\:\Gamma\to Y$ 
that is metric-minimizing relative to $A$ and $f\succcurlyeq h\rel A$.
\end{thm}

\parit{Proof.}
Let us parametrize each edge of $\Gamma$ by $[0,1]$.
A map $h\:\Gamma\to Y$ will be called \emph{straight} if it
sends each edge of $\Gamma$ to a constant-speed geodesic path in $Y$.

If $h\:\Gamma\to Y$ is straight, then $f\succcurlyeq h\rel A$ if and only if 
\[|f(v)-f(w)|_Y\ge |h(v)-h(w)|_Y\]
for any two adjacent vertexes $v$ and $w$ in $\Gamma$.
In particular, we can assume that the given map $f$ is straight.

By finiteness of the number of vertexes and Zorn's lemma,
it is sufficient to prove that for any ordered sequence of straight maps $f_1\succcurlyeq f_2\succcurlyeq \dots$ there exists a map $f\preccurlyeq f_n$ for all $n$.

Assume the contrary; let us apply the \emph{ultralimit+projection construction}.

Namely, fix an ultrafilter $\omega$; denote by $Y^\omega$ the ultrapower of $Y$. 
Then $Y^\omega$
is a $\CAT(0)$ space that contains $Y$ as a closed convex subset. 
The $\omega$-limit $f_\omega\:\Gamma\to Y^\omega$ is well-defined since all $f_n$ are Lipschitz continuous. 
Denote by $f'$ the composition of $f_\omega$ with the nearest point projection $Y^\omega\to Y$.

The nearest point projection to a closed convex set in $\CAT(0)$ space is short.
Therefore $f_n\succcurlyeq f'$ for any $n$.
Let $f''$ denote the straightening of $f'$.
Then $f'\succcurlyeq f''$ and the claim follows.
\qeds

\begin{thm}{Proposition}\label{prop:metric-min-graph}
Let $Y$ be a $\CAT(0)$ space, 
$\Gamma$ be a finite  graph and $A$ be a subset of its vertexes.
Assume that the assignment $v\mapsto v'$ is a metric-minimizing map $\Gamma\to Y$ relative to $A$.
Then
\begin{enumerate}[(a)]
\item each edge of $\Gamma$ maps to a geodesic;
\item\label{prop:metric-min-graph:x} for any vertex $v\notin A$ and any $x\ne v'$
there is an edge  $[v,w]$ in $\Gamma$ such that
$\measuredangle[v'\,^{w'}_x]\ge \tfrac\pi2$;
\item\label{sum>=2pi} for any vertex $v\notin A$ and any cyclic order $w_1,\dots,w_n$ of adjacent vertexes we have
\[\measuredangle[v'\,^{w'_1}_{w'_2}]+\dots+\measuredangle[v'\,^{w'_{n-1}}_{w'_n}]+\measuredangle[v'\,^{w'_n}_{w'_1}]\ge 2\cdot\pi.\]
\end{enumerate}
\end{thm}

\begin{wrapfigure}{r}{22 mm}
\begin{lpic}[t(-0 mm),b(-0 mm),r(0 mm),l(0 mm)]{pics/not-sufficient(1)}
\end{lpic}
\end{wrapfigure}

\parit{Remark.}
The conditions in the proposition do not guarantee that the map $f$ is metric-minimizing.
An example can be guessed from the diagram, where the solid points form the set~$A$.

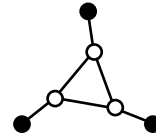
\parit{Proof.}
The first condition is evident.

Assume the second condition does not hold at a vertex $v\zz\notin A$;
that is, there is a point $x\in Y$ such that
$\measuredangle[v'\,^{w'}_x]< \tfrac\pi2$
for any adjacent vertex $w$.
In this case, moving $v'$ toward $x$ along $[v',x]$ decreases the lengths of all edges adjacent to $v$, a contradiction.

Assume the third condition does not hold; 
that is, the sum of the angles around a fixed interior vertex $v'$ is less than $2\cdot\pi$.

Recall that the space of directions $\Sigma_{v'}$ is a $\CAT(1)$ space.
Denote by $\xi_1,\dots,\xi_n$ the directions of $[v',w'_1],\dots, [v',w'_n]$ in $\Sigma_{v'}$.
By assumption, we have
\[|\xi_1-\xi_2|_{\Sigma_{v'}}+\dots+|\xi_n-\xi_1|_{\Sigma_{v'}}<2\cdot\pi.\]
By Reshetnyak's majorization theorem,
the closed broken line $[\xi_1,\dots,\xi_k]$ is majorized by a convex spherical polygon $P$.

Note that $P$ lies in an open hemisphere with a pole at some point in $P$.
Choose $x\in Y$ so that the direction from $v'$ to $x$ coincides with the image of the pole in $\Sigma_{f(v)}$.
This choice of $x$ contradicts (\emph{\ref{prop:metric-min-graph:x}}).
\qeds

\section{Key lemma}\label{Key Lemma}

\begin{thm}{Lemma}\label{lem:graph}
Let $Y$ be a $\CAT(0)$ space and $s\:\DD\to Y$ 
be a metric-minimizing disc.
Assume $F\subset \DD$ is a finite set such that $\hat\pi_s(F)$ has finite
diameter in $\<\DD\>_s$.
Then there exists a finite piecewise geodesic graph $\Gamma$ embedded in $\<\DD\>_s$ that contains a geodesic between any pair of points in $\hat\pi_s(F)$.
\end{thm} 

\parit{Proof.} 
For any pair $x,y\in F$, connect $\hat x$ to $\hat y$ by a minimizing geodesic in $\<\DD\>_s$. 
We can assume that the constructed geodesics 
are either disjoint or their intersection is formed by finite collections of arcs and points.

Indeed, if some number of geodesics $\gamma_1,\dots,\gamma_n$ already has this property and we are given points $x$ and $y$, then
we choose a minimizing geodesic $\gamma_{n+1}$ from $x$ to $y$ that maximizes the time it spends in $\gamma_1,\dots,\gamma_n$  
in the order of importance.
Namely, 
\begin{itemize}
\item  among all minimizing geodesics connecting $x$ to $y$
choose one that spends maximal time in $\gamma_1$ --- in this case, $\gamma_{n+1}$ intersects $\gamma_1$ along the empty set, 
a one-point set, or a closed arc.
\item among all minimizing geodesics as above
choose one that spends maximal time in $\gamma_2$ --- in this case, $\gamma_{n+1}$ intersects $\gamma_2$ along at most two arcs and points.
\item and so on.
\end{itemize}

It follows that together the constructed geodesics form a finite graph $\Gamma$ as required.
\qeds

\begin{thm}{Key lemma}\label{lem:key}
Let $Y$ be a $\CAT(0)$ space and $s\:\DD\to Y$ 
be a metric-minimizing disc.
Given a finite set $F\subset \DD$
there is 
\begin{enumerate}[(1)]
	\item a $\CAT(0)$ disc-retract $W$ with boundary curve $\delta$;
	\item a map $p\:F\to W$ such that
\[|p(x)-p(y)|_W\le \<x-y\>_s\] 
for $x,y\in F$ and $p(x)=\delta(x)$ for $x\in F\cap \partial\DD$;
  \item a short map $q\:W\to Y$ such that
\[s(x)=q\circ p(x)\] 
for any $x\in\partial\DD\cap F$.
\end{enumerate}
 
\end{thm} 

\parit{Proof.} If $\partial \DD\cap F= \emptyset$,
then one can take a one-point space as $W$ and arbitrary maps $p\:F\to W$ and $q\:W\to Y$.
So suppose $\partial \DD\cap F\ne\emptyset$.

Without loss of generality, we may assume that the distance $\<x-y\>_s$
between any pair of points $x,y\in F$ is finite.
Indeed, since the boundary curve $s|_{\partial\DD}$ is rectifiable,
this always holds for pairs of points in $\partial \DD\cap F$.
Consider the subset $F'\subset F$ that lies at finite $\<{*}-{*}\>_s$-distance from one (and therefore any) point in $\partial \DD\cap F$.
Suppose $p'\:F'\to W$ and $q\:W\to Y$ are maps satisfying the proposition for $F'$.
Extend $p'$ to $F$ by sending $F\backslash F'$ to one point in $W$.
The resulting map $p$ together with $q$ will then satisfy the proposition for $F$.

By Lemma~\ref{lem:graph}, there exists a finite piecewise geodesic graph $\Gamma$ embedded in $\<\DD\>_s$ that contains $F$ as a subset of its vertexes.
According to Proposition~\ref{prop:|D|},
 $\tau_s$ embeds $\Gamma$ in $|\DD|_s$.
By Proposition~\ref{prop:|D|},
$|\DD|_s$ is a disc-retract.
Therefore $\Gamma$ can be (and will be) considered as a graph embedded into the plane.

By Proposition~\ref{prop:metric-min-graph-exist}, there is a map 
$u\:\Gamma\to Y$ metric-minimizing relative to $A=F\cap\partial\DD$ such that
\[s|_\Gamma\succcurlyeq u\rel A.\eqlbl{eq:>=}\]

Fix an open disc $\Delta$ cut by $\Gamma$ from $|\DD|_s$.
By Reshetnyak's theorem, the closed curve $u|_{\partial\Delta}$
is majorized by a convex plane polygon, possibly degenerating to a point or a line segment.
Note that the angle of the majorizing polygon cannot be smaller than the angle between the corresponding edges in $u(\Gamma)\subset Y$.

Let us glue the majorizing polygons into $\<\Gamma\>_u$;
denote by $W$ the resulting space.
According to Proposition~\ref{prop:metric-min-graph}(\ref{sum>=2pi}), the angle around each inner vertex has to be at least $2\cdot\pi$.
Clearly, $W$ is a disc-retract;
in particular, it is simply connected.
It follows that $W$ is a $\CAT(0)$ space.

The short map $q\:W\to Y$ is constructed by gluing together the maps provided by Reshetnyak's majorization theorem.
The space $W$ comes with a natural short map $\<\Gamma\>_u\to W$.

Define $p(x)$ for $x\in F$ as the image of the corresponding vertex of $\Gamma$ in $W$.
By \ref{eq:>=}, 
\[|p(x)-p(y)|_W\le \<x-y\>_s\]
for any $x,y\in F$.

By construction, the pair of maps $p,q$ meet all conditions.
\qeds

The following establishes a connection between the 
key lemma (\ref{lem:key}) and the extension lemma (\ref{lem:finite-whole}).

\begin{thm}{Lemma}\label{lem:S-closed}
Let $Y$ be a $\CAT(0)$ space and $s\:\DD\to Y$ 
be a metric-minimizing disc. Let $W$ be a
$\CAT(0)$ disc-retract with boundary curve $\delta$. For a given 
finite set $F\subset \DD$ we define $\mathfrak{S}_F$ to be the family of maps  
 $p\:F\to W$ such that
\[|p(x)-p(y)|_W\le \<x-y\>_s\] 
for $x,y\in F$ with $p(x)=\delta(x)$ for $x\in F\cap \partial\DD$ and such that there exists
a short map $q\:W\to Y$ with
\[s(x)=q\circ p(x)\] 
for any $x\in\partial\DD\cap F$.
Then $\mathfrak{S}_F$ is closed under pointwise convergence. 
\end{thm}

The proof is an application of the ultralimit+projection construction;
we used it before and will use it again later.

\parit{Proof.}
Consider a converging sequence $p_n\in  \mathfrak{S}_F$;
denote by $p_\infty$ its limit.
For each $p_n$ there is a short map $q_n\:W\to Y$ satisfying the condition above.
Pass to its ultralimit $q_\omega\:W\to Y^\omega$.
Recall that $Y$ is a closed convex set in $Y^\omega$.
In particular, the nearest point projection $\nu\:Y^\omega\to Y$ is well-defined and short.
Therefore, the composition $q=\nu\circ q_\omega$ is short.
Finally note that the maps $p_\infty\:F\to W$ and $q\:W\to Y$ satisfy the condition above.
\qeds

\section{Extension lemma}\label{Finite-whole extension lemma}

\begin{thm}{Extension lemma}\label{lem:finite-whole}
Suppose that $X$ is a set 
and $Y$ is a compact topological space.
Assume that for any finite set $F\subset X$ 
a nonempty set $\mathfrak{S}_F$ of maps  $F\to Y$ is given, such that
\begin{itemize}
\item $\mathfrak{S}_F$ is closed under pointwise convergence;
\item for any subset $F'\subset F$ and any map $h\in \mathfrak{S}_F$
the restriction $h|_{F'}$ belongs to $\mathfrak{S}_{F'}$. 
\end{itemize}

Then there is a map $h\: X\to Y$ such that $h|_F\in \mathfrak{S}_F$ for any finite set $F\subset X$.
\end{thm}

\parit{Proof.}
Consider the space $Y^X$ of all maps $X\to Y$ equipped with the product topology.

Denote by $\bar{\mathfrak{S}}_F$ the set of maps $h\in Y^X$ such that its restriction $h|_F$ belongs to $\mathfrak{S}_F$.
By assumption, the sets $\bar{\mathfrak{S}}_F\subset Y^X$ are closed and any finite intersection of these sets is nonempty.

According to Tikhonov's theorem, $Y^X$ is compact.
By the finite intersection property, the intersection $\bigcap_F\bar{\mathfrak{S}}_F$ for all finite sets $F\subset X$ is nonempty.
Hence the statement follows.
\qeds

Note that if $X$ and $Y$ are metric spaces and $A$ is a subset in $X$,
then one can take as $\mathfrak{S}_F$ the short maps $F\to Y$ that coincide with a given short map $A\to Y$ on $A\cap F$.
This way we obtain the following corollary; it is closely related to \cite[Proposition 5.2]{lang-shroeder}.
In a similar fashion, we will use the lemma in the proof of our main theorem.

\begin{thm}{Corollary}
Let $X$ and $Y$ be metric spaces, $A\subset X$, and $f\:A\to Y$ a short map.
Assume $Y$ is compact and for any finite set $F\subset X$ there is a short map $F\to Y$ that agrees with $f$ in $F\cap A$.
Then there is a short map $X\to Y$ that agrees with $f$ in $A$.
\end{thm}

\section{Proof assembling}\label{Main theorem}

\parit{Proof of the main theorem.}
Given a finite set $F\subset \DD$,
denote by $\mathcal{W}_F$
the set of isometry classes of spaces $W$ that meet the conditions of the key lemma~(\ref{lem:key})
for $F$;
according to this lemma, $\mathcal{W}_F\ne\emptyset$.
Note that for two finite sets $F\subset F'$ in $\DD$,
we have $\mathcal{W}_F\supset \mathcal{W}_{F'}$.

According to the compactness lemma (\ref{lem:compact}), $\mathcal{W}_F$ is compact.
Therefore 
\[\mathcal{W}
=
\bigcap_{F}\mathcal{W}_F\ne \emptyset\]
where the intersection is taken over all finite subsets $F$ in $\DD$.

Fix a space $W$ from $\mathcal{W}$;
the space $W$ is a $\CAT(0)$ disc-retract,
such that given a finite set $F\subset \DD$ there is a map $h_F\:F\to W$ that is short with 
respect to $\<{*}-{*}\>_s$ 
and a short map $q_F\:W\to Y$ such that $q_F\circ h_F$ agrees with $s$ on $\partial\DD\cap F$.

Given a finite set $F\subset \DD$,
denote by $\mathfrak{S}_F$ the set of all maps $h_F\:F\to W$ described above.

By Lemma \ref{lem:S-closed}, $\mathfrak{S}_F$ is closed.
The condition on the restriction of $h_F\in  \mathfrak{S}_F$ in the extension lemma (\ref{lem:finite-whole}) is evident.
Applying the lemma,
we get a map $h\:\DD\to W$ such that $h|_F\in \mathfrak{S}_F$
for any finite set $F\subset \DD$.

Our next aim is to show that there is a single map $q$ such that
for all finite sets $F$ the composition $q\circ h|_F$ agrees with
$s$ on $\partial\DD\cap F$.
This is done by applying the ultralimit+projection construction:

Choose a sequence of finite sets $F_n$ such that $F_n$ get denser and denser in $\DD$ and
the intersections $F_n\cap\partial \DD$ get denser and denser in $\partial \DD$; 
denote by $q_n$ the corresponding maps.
Let $q_\omega\:W\to Y^\omega$ be the ultralimit of $q_n$ and set $q=\nu\circ q_\omega$,
where $\nu\:Y^\omega\to Y$ is the nearest point projection.
By construction, $q\:W\to Y$ is short and $q\circ h$ agrees with $s$ on $\partial \DD$.
Note that we cannot conclude $s\succcurlyeq q\circ h$ because $h$ might not be continuous.

By construction, the map $h$ induces a short map $\hat h\:\<\DD\>_s\to W$ 
such that $\hat h\circ\delta_s$ is the boundary curve of $W$.
By Lemma \ref{lem:maj is isom}, $\hat h$ is an isometry and the statement follows.
\qeds

\section{About the old theorem}\label{sec:old-thm}

In this section we formulate and prove two versions of the main theorem in \cite{petrunin-metric-min} with corrections as in \cite{petrunin-metric-min-correction}.

Note that the meaning of the term \emph{metric-minimizing} in the present paper differs from its meaning in \cite{petrunin-metric-min} --- the old paper used a weaker definition and the theorem requires an additional assumption.

In Section~\ref{sec:defs}, we introduced three ways to pull back the metric along a map $f\:X\to Y$ from a topological space $X$
 to a metric space $Y$: 
the induced length metric
$\langle x-y\rangle_f$,
the induced intrinsic metric $\langle| x-y|\rangle_f$
and the induced connecting metric $|x-y|_f$.
From the definitions, we have that 
\[\langle x-y\rangle_f\ge \langle| x-y|\rangle_f\ge |x-y|_f\]
for any $x,y\in X$.
The second inequality is strict for generic maps, 
an example of a map $f$ with strict first inequality is given in \cite[4.2]{petrunin-intrinisic}.
(If $f$ is an embedding, then equality holds \cite[4.5]{ledonne}.)
\begin{figure}[h]
\centering
\begin{tikzpicture}[scale=1.5]

\node (0) at (1,.5) {$\langle|X|\rangle_f$};
  \node (1) at (1,-.5) {$|X|_f$};
  \node (2) at (1,1.5){$\<X\>_f$};
  \node (11) at (3,1){$Y$};
  \node (12) at (-1,1) {$X$};
\draw[
    >=latex,
    auto=right,                      
    loop above/.style={out=75,in=105,loop},
    every loop,
    ]
   (2) edge node{$\hat\tau_f$}(0)
   (0) edge node{$\bar\tau_f$}(1)
   (12) edge[bend right=10] node[swap]{$\hat{\bar \pi}_f$}(0)
   (0) edge[bend right=10] node[swap]{$\hat{\bar f}$}(11)
   (12) edge[bend left=90] node[swap]{$f$}(11)
   (12) edge[bend right=-10] node[swap]{$\hat\pi_f$}(2)
   (2) edge[bend right=-10] node[swap]{$ \hat f$}(11)
   (12) edge[bend right] node{$\bar \pi_f$}(1)
   (1) edge[bend right] node{$\bar f$}(11);
\end{tikzpicture}
\end{figure}
The diagram is an extension of the diagram on page~\pageref{diagram-page} which includes $\langle|X|\rangle_f$.
From above, both maps $\hat\tau_f$ and $\bar\tau_f$ are short, and $\tau_f=\bar\tau_f\circ\hat\tau_f$.
Note that $\tau_f$ might fail to be injective while $\bar\tau_f$ is always injective.

Our first formulation uses the $\langle|{*}-{*}|\rangle$-metric instead of $\langle{*}-{*}\rangle$ which was used in the old formulation.

\begin{thm}{Conjecture for $\bm{\langle|{*}-{*}|\rangle}$}\label{thm:old1}
Let $Y$ be a $\CAT(0)$ space and $s\:\DD\to Y$ a continuous map that satisfies the following property: 
if a continuous map $s'\:\DD\to Y$ agrees with $s$ on $\partial\DD$ and
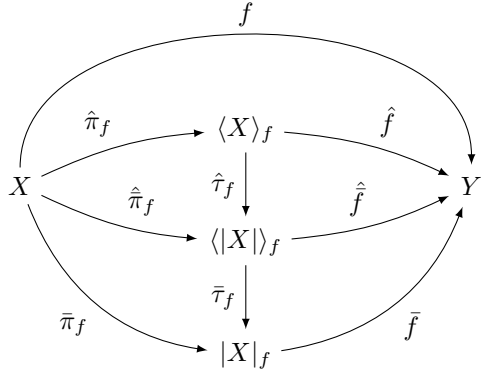
\[\<|x-y|\>_{s'}\le \<|x-y|\>_{s}\]
for any $x,y\in \DD$,
then equality holds for all pairs of $x$ and $y$.
Assume that the space $\<|\DD|\>_s$ is compact.
Then $\<|\DD|\>_s$ is $\CAT(0)$.
\end{thm}

In the published paper, it was formulated as a theorem.
However, as pointed out by Alexander Lytchak, our proof contained a gap.
At least after all these years, we were not able to reconstruct the proof from what was written there.

The second formulation uses the metric $\langle{*}-{*}\rangle$ as in the original formulation.

\begin{thm}{Old theorem for $\bm{\langle{*}-{*}\rangle}$}\label{thm:old2}
Let $Y$ be a $\CAT(0)$ space and $s\:\DD\to Y$ a continuous map that satisfies the following property: 
if a continuous map $s'\:\DD\to Y$ agrees with $s$ on $\partial\DD$ and
\[\<x-y\>_{s'}\le \<x-y\>_{s}\]
for any $x,y\in \DD$,
then equality holds for all pairs of $x$ and $y$.
Assume that the function $(x,y)\mapsto\<x-y\>_{s}$ is continuous.
Then $\<\DD\>_s$ is $\CAT(0)$.
\end{thm}

Note that continuity of the function $(x,y)\mapsto\<x-y\>_{s}$
implies that $\<\DD\>_s$ is compact.
Therefore the former condition is stronger than the latter.  
The sketch of proof given in \cite{petrunin-metric-min} implicitly used that the metric $(x,y)\zz\mapsto\<x-y\>_{s}$ is continuous.
(We do not know if compactness of $\<\DD\>_s$ alone is sufficient.)

Theorem~\ref{thm:old2} follows from the following proposition and Theorem~\ref{thm:old1}.
 
\begin{thm}{Proposition}\label{prop:<||>=<>}
Let $Y$ be a metric space and let $s\:\DD\to Y$ be a continuous map without bubbles.
Assume that the function $(x,y)\mapsto \<x-y\>_s$ is continuous.
Then 
\[\<x-y\>_s=\<|x-y|\>_s\]
for any $x,y\in \DD$.
\end{thm}

Before going into the proof, let us give an example showing that the proposition is not trivial.

Consider a pseudo-arc $P\subset \DD$ and let $s$ be the quotient map $\DD\to \DD/P$.
Evidently, 
\[\<|x-y|\>_s=0\]
for any $x,y\in P$.
However, since $P$ contains no curves it is not at all evident that 
\[\<x-y\>_s=0\]
for any $x,y\in P$.

We present an argument of Taras Banakh \cite{banakh};
it works in dimension two but we are not aware of a generalization to higher dimensions.

\begin{thm}{Lemma}\label{lem:subdivision}
Let $Y$ be a metric space and $s\:\DD\to Y$ a continuous map without bubbles.
Assume that the function $(x,y)\mapsto \<x-y\>_s$ is continuous.
Then there is a finite collection of curves in $\DD$ with finite total $\<{*}-{*}\>$-length 
that divide $\DD$ into subsets with arbitrarily small $\<{*}-{*}\>$-diameter.
\end{thm}

\parit{Proof.}
Note that given $\eps>0$ there is $\delta>0$ such that a set of diameter $\delta$ in $\<\DD\>_s$ can not separate a set of diameter at least $\eps$ from the boundary curve.
If this is not the case, 
then sets of arbitrarily small diameter can separate a set of diameter at least $\eps$ from the boundary.
Passing to a limit we get a one-point set that separates a set from the boundary curve; that is, $s$ has a bubble --- a contradiction. 

Let us subdivide $\DD$ into small pieces by curves, say by vertical and horizontal lines.
Since the function $(x,y)\mapsto \<x-y\>_s$ is continuous,
we can assume that all pieces have small $\<{*}-{*}\>_s$-diameter;
that is, given $\eps>0$ we can assume that $\<x-y\>_s<\eps$ for any two points $x$ and $y$ in one piece.

It remains to modify the decomposition to make the boundary curves $\<{*}-{*}\>_s$-rectifiable.

Subdivide the curves into arcs with $\<{*}-{*}\>_s$-diameter smaller than $\tfrac{\delta}{5}$ and exchange this piece with a curve of $\<{*}-{*}\>_s$-length smaller than $\tfrac{\delta}{5}$.
The new arc together with the old one form a set of  $\<{*}-{*}\>_s$-diameter at most $\delta$.\footnote{Note that this arc might travel far in the Euclidean metric on $\DD$.}
Therefore we might add to a piece a subset of diameter at most $\eps$ and the total $\<{*}-{*}\>_s$-diameter of each piece remains below $3\cdot\eps$.
The curves might cut more pieces from $\DD$, but by the same argument each of these pieces will have $\<{*}-{*}\>_s$-diameter below $3\cdot\eps$.
\qeds

\parit{Proof of \ref{prop:<||>=<>}.}
Fix points $x,y\in \DD$.
It is sufficient to construct a path $\alpha$ from $x$ to $y$ in $\DD$ such that the length of $s\circ\alpha$ is arbitrarily close to $\<|x-y|\>_s$.

Since $(x,y)\mapsto \<x-y\>_s$ is continuous, $\<\DD\>_s$ and therefore $\<|\DD|\>_s$ are compact.
In particular, there is a minimizing geodesic $\gamma$ from $\hat{\bar x}=\hat{\bar \pi}(x)$ to $\hat{\bar y}=\hat{\bar \pi}(y)$ in $\<|\DD|\>_s$.
Denote by $\Gamma$ the inverse image of $\gamma$ in $\DD$;
this is a connected compact set that does not have to be path-connected.

To construct the needed path $\alpha$, it is sufficient to prove the following claim:

\begin{itemize}
 \item[$\bigstar$] Given $\eps>0$ there is a set $\Gamma'\subset \Gamma$ 
 and a collection of paths $\alpha_0,\dots,\alpha_n$ such that 
 (1) the total length of $s(\alpha_i\backslash\Gamma)$ is at most $\eps$, 
 (2) The set $\Gamma'$ is a union of a finite collection of closed connected sets $\Gamma_0,\dots,\Gamma_n$, 
 (3) diameter of each $\Gamma_i$ is at most $\eps$, 
 (4) $x\in\alpha_0$, $y\in\alpha_n$, and 
 (5) the union $\Gamma'\cup \alpha_0\cup\dots\cup\alpha_n$ is connected.
\end{itemize}

Indeed, once Claim $\bigstar$ is proved, one can apply it recursively for a sequence of $\eps_n$ that converges to zero very fast.
Namely, we can apply the claim to each of the subsets $\Gamma_n$ and take as $\Gamma''$ the union of all closed subsets provided by the claim.
This way we obtain a nested sequence of closed sets $\Gamma\supset \Gamma'\supset\Gamma''\supset\dots$ which break into a finite union of closed connected subsets of arbitrarily small diameter and a countable collection of arcs with total length at most $\eps_1+\eps_2+\dots$ 
outside of $\Gamma$.
Set 
\[\Phi=\Gamma\cap \Gamma'\cap\Gamma''\cap\dots\]
Note that there is a simple curve from $x$ to $y$ that runs in the constructed arcs and $\Phi$.
The part of the curve in $\Gamma$ contributes at most $\<|x-y|\>_s$ to its $\<{*}-{*}\>_s$-length.
Therefore the total length of the curve can not exceed 
\[\<|x-y|\>_s+\eps_1+\eps_2+\dots;\]
hence the result will follow.

It remains to prove $\bigstar$.

Fix a subdivision $\Upsilon_1,\dots,\Upsilon_k$ of $\DD$ provided by Lemma~\ref{lem:subdivision} for the given $\eps$.
Denote by $\Delta$ the union of all the cutting curves.

By the regularity of $\<{*}-{*}\>_s$-length, we may cover $\Delta\cap\Gamma$ by a finite collection of arcs with total  $\<{*}-{*}\>_s$-length arbitrarily close to $\<{*}-{*}\>_s$-length of $\Delta\cap\Gamma$.
Denote these arcs by $\alpha_0,\dots,\alpha_n$.
Without loss of generality, we may assume that $x\in\alpha_0$ and $y\in\alpha_n$.

Consider a finite graph with vertexes labeled by $\alpha_0,\dots,\alpha_n$;
two vertexes $\alpha_i$ and $\alpha_j$ are connected by an edge if there is a connected set $\Theta\subset \Gamma\cap\Upsilon_k$ 
for some $k$ such that $\Theta$ intersects $\alpha_i$ and $\alpha_j$.
Note that the graph is connected. Therefore we may choose a path from $\alpha_0$ to $\alpha_n$ in the graph.

The path corresponds to a sequence of arcs $\alpha_i$ and a sequence of $\Theta$-sets.
The $\Theta$-sets that correspond to the edges in the path can be taken as $\Gamma_i$ in the claim.
Hence the claim and therefore the proposition follow. 
\qeds

\section{Saddle surfaces}\label{sec:smooth}

In this section we will discuss the relation between metric-minimizing discs and saddle discs.

Recall that a map $s\:\DD\to \RR^m$ is called \emph{saddle} if for any hyperplane $\Pi\subset\RR^m$ each of the connected components of $\DD\backslash s^{-1}\Pi$ meets the boundary.

If $s$ is a smooth embedding in $\RR^3$, then it is saddle if and only if the  obtained surface has nonpositive Gauss curvature. 
An old conjecture of Samuel Shefel states that any saddle disc in a $\RR^3$ is $\CAT(0)$ with respect to its length metric, see \cite{shefel-3D}.

It is evident that any metric-minimizing disc $s$ in a Euclidean space is saddle.

\parbf{Three-dimensional case.}
In general, a saddle disc may not be globally metric-minimizing.
An example is shown in the picture.
It is a saddle polyhedral disc made from 10 triangles
with a hexagon boundary.
The boundary curve goes along the Y-shape marked with bold lines;
each segment of the Y-shape is traveled twice back and forth.

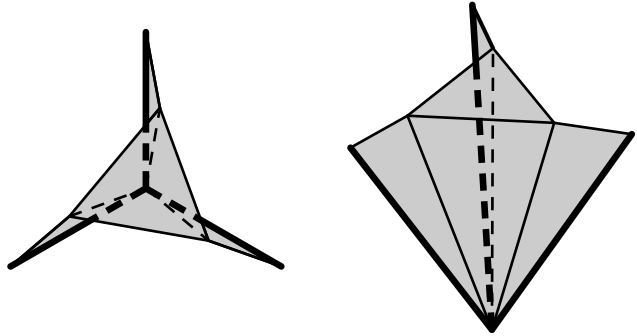
\begin{figure}[ht!]
\centering
\begin{lpic}[t(-0 mm),b(-0 mm),r(0 mm),l(0 mm)]{pics/not-sufficient-disc(1)}
\end{lpic}
\end{figure}

The picture is rotationally symmetric by the angle $\tfrac23\cdot\pi$.
A shortening deformation can be obtained by rotating the central triangle slightly counterclockwise and extending the map on the remaining 9 triangles linearly.

By smoothing this example one can produce a smooth saddle disc that is not metric-minimizing.

Proposition~\ref{prop:smooth} states that there are no local examples of that type that are
smooth and \emph{strictly saddle} meaning that the principal curvatures at each interior point have opposite signs.
To prove this proposition we need to introduce a certain energy functional.

Let $s\:\DD\to\RR^3$ be a smooth map.

Fix an array of vector fields $\bm{v}=(v_1,\dots,v_k)$ on $\DD$. 
Assume that each integral curve of vector fields $v_i$ goes from boundary to boundary yielding sweep outs of the whole disc $\DD$. 

Consider the energy functional 
\[E_{\bm{v}}s
:=
\sum_i\int\limits_\DD |v_is|^2,\]
where $vs=ds(v)$ denotes the derivative of $s$ in the direction of the field $v$.
Set 
\[\Delta_{\bm{v}}s=\sum_iv_i(v_is).\]
It is convenient to think of the operator $s\mapsto \Delta_{\bm{v}}s$
as an analog of the Laplacian.

Note that 
\begin{enumerate}[(i)]

\item $E_{\bm{v}}$ is well-defined for any Lipschitz map $s$.

\item $E_{\bm{v}}$ is convex; that is,
\[E_{\bm{v}}s_t
\le 
(1-t)\cdot E_{\bm{v}} s_0+t\cdot E_{\bm{v}} s_1,\]
where $s_t=(1-t)\cdot s_0+t\cdot s_1$ and $0\le t\le 1$.
Moreover, the equality holds for any $t$ if and only if for any $i$ we have $v_is_0=v_is_1$ almost everywhere.

\item
If $s_0$ is a smooth $E_{\bm{v}}$-minimizing map in the class of Lipschitz maps with given boundary data, then $s_0$ is metric-minimizing. 
\\
Indeed, if $s_1\preccurlyeq s_0$, then $s_1$ has to be Lipschitz.
It follows that $E_{\bm{v}} s_1\le E_{\bm{v}} s_0$ and from convexity $E_{\bm{v}} s_t\le E_{\bm{v}} s_0$ if $0\le t\le 1$.
Since $s_0$ is $E_{\bm{v}}$-minimizing, $E_{\bm{v}} s_t= E_{\bm{v}} s_0$ for any $t$.
Hence $v_is_0=v_is_1$ almost everywhere.
Since $\DD$ is swept out by arcs of integral curves of $v_i$, the latter implies $s_0=s_1$.

\item A smooth map $s\:\DD\to\RR^3$ is a $E_{\bm{v}}$-minimizing map among the class of Lipschitz maps with given boundary if and only if
\[\Delta_{\bm{v}}s=0.\]

\end{enumerate}

The discussion above reduces Proposition~\ref{prop:smooth} to the following.

\begin{thm}{Claim}
Assume $s\:\DD\to \RR^3$ is a smooth strictly saddle surface. 
Then for any interior point $p\in\DD$ there is an array of 4 vector fields $\bm{v}=(v_1,v_2,v_3,v_4)$ such that the equation \[\Delta_{\bm{v}}s=0\eqlbl{eq:laplasian}\]
holds in an open neighborhood of $p$.
\end{thm}

\parit{Proof.}
Denote 
by $\kappa_1,\kappa_2$ the principal curvatures,
and by $e_1,e_2$ the corresponding unit principal vectors. 
Further, denote by $a_1,a_2$ a pair of asymptotic vectors; that is, the normal curvatures in these directions vanish. 
We can assume that $a_1,a_2$ form coordinate vector fields in a neighborhood of $p$.

Set $v_1=\tfrac 1{\sqrt{|\kappa_1|}}\cdot e_1$ and $v_2=\tfrac 1{\sqrt{|\kappa_2|}}\cdot e_2$. 
It remains to show that one can choose smooth functions  $\lambda_1$ and $\lambda_2$ 
so that \ref{eq:laplasian}
holds in a neighborhood of $p$ for $v_3=\lambda_1\cdot a_1$ and $v_4=\lambda_1\cdot a_1$.

Note that the sum $v_1(v_1s)+v_2(v_2s)$ has vanishing normal part.
That is \[v_1(v_1s)+v_2(v_2s)\] is a tangent vector to the surface.

Since $a_i$ are asymptotic,
the vectors $a_1(a_1s)$ and $a_2(a_2s)$ have vanishing normal parts.
Therefore, for any choice of $\lambda_i$,
the following two vectors are also tangent
\begin{align*}
v_3(v_3s)&=\lambda_1^2\cdot a_1(a_1s)+\tfrac12\cdot a_1\lambda_1^2\cdot a_1s
\\
v_4(v_4s)&=\lambda_2^2\cdot a_2(a_2s)+\tfrac12\cdot a_2\lambda_2^2\cdot a_2s.
\end{align*}

Set $w=(\lambda_1^2,\lambda_2^2)$.
Note that the system \ref{eq:laplasian} can be rewritten as 
\[\left(\begin{smallmatrix}
   1&0\\0&0
  \end{smallmatrix}\right)
w_x
+
\left(\begin{smallmatrix}
   0&0\\0&1
  \end{smallmatrix}\right)
w_y=h(x,y,w),\]
where $h\:\RR^3\to\RR^2$ is a smooth function.

Change coordinates by setting $x=t+z$ and $y=t-z$.
Then the system takes the form 
\[w_t+\left(\begin{smallmatrix}
   1&0\\0&-1
  \end{smallmatrix}\right)
w_z=h(t+z,t-z,w),\]
which is a semilinear hyperbolic system.
According to \cite[Theorem 3.6]{bressan}, it can be solved locally for smooth initial data at $t=0$.

It remains to choose $v_3$ and $v_4$ for solution such that $\lambda_1, \lambda_2>0$ in a small neighborhood of $p$.
\qeds

\parbf{Four-dimensional case.}
Except for constructing an energy as we did above,
we do not see any way to show that a given smooth surface is metric-minimizing.
Locally, the appropriate energy functional can be described by three functions defined on the disc.
These three functions are subject to certain differential equations.
Straightforward computations show that on generic smooth saddle surfaces in $\RR^4$ 
there is no solution even locally.

For that reason we expect that generic smooth saddle surfaces in $\RR^4$ are not locally metric-minimizing. 
That is, arbitrarily small neighborhoods of any point admit deformations that shrink 
the length metric and keep the boundary fixed.
On the other hand, we do not have an example of a saddle surface for which this condition would hold at a single point.

\Addresses


\begin{thebibliography}{52}

\bibitem{A}
\begin{otherlanguage}{russian}
Александров, А. Д., 
\textit{Линейчатые поверхности в метрических пространствах.}
Вестник ЛГУ 2 (1957): 15---44.
\end{otherlanguage}


\bibitem{AB} Alexander, S.; Bishop, R., \textit{Gauss equation and injectivity radii for subspaces in spaces of curvature bounded above.} 
Geom. Dedicata 117 (2006), 65--84. 


\bibitem{akp}
Alexander, S.; Kapovitch, V. and Petrunin, A.,
\textit{An invitation to Alexandrov geometry: CAT(0) spaces.}
SpringerBriefs in Mathematics, 2019. 

\bibitem{banakh} Banakh, T. \textit{Running most of the time in a connected set},  MathOverflow   \texttt{https://mathoverflow.net/q/308172}, 

\bibitem{bressan} Bressan, A.,
\textit{Hyperbolic systems of conservation laws.
The one-dimensional Cauchy problem.}
Oxford Lecture Series in Mathematics and its Applications, 20, 2000.

\bibitem{BBI}Burago, D.; Burago, Y. and Ivanov, S.,
\textit{A course in metric geometry.}
Graduate Studies in Mathematics, 33, 2001.

\bibitem{daverman} Daverman, R. J.,
\textit{Decompositions of Manifolds}
Academic Press Inc., N.Y., 1986.

\bibitem{G} Gromov, M., \textit{Singularities, expanders and topology of maps. Part 1 : Homology versus volume in the spaces of cycles.} Geom. Funct. Anal. 19 (2009) no. 3, 743--841.

\bibitem{HK}
Hu, T.; Kirk, W. A.,  \textit{Local contractions in metric spaces.}
Proc. Amer. Math. Soc. 68.1 (1978), 121--124.


\bibitem{KF}
\begin{otherlanguage}{russian}
Колмогоров, А. Н.;
Фомин, С. В.,
\textit{Элементы теории функций и функционального анализа.}
Издание седьмое, 2004.
\end{otherlanguage}

\bibitem{lang-shroeder} Lang, U.; Schroeder, V.
\textit{Kirszbraun's theorem and metric spaces of bounded curvature.}
Geom. Funct. Anal. 7 (1997), no. 3, 535--560. 

\bibitem{ledonne} Le Donne, E.
\textit{Lipschitz and path isometric embeddings of metric spaces.} 
Geom. Dedicata 166 (2013), 47--66.

\bibitem{KS}Korevaar, N. J.; Schoen, R. M. ``Sobolev spaces and harmonic maps for metric space targets,'' Comm. Anal. Geom., 1(3-4):561--659, 1993.

\bibitem{L} Lytchak, A, 
\textit{On the geometry of subsets of positive reach.} 
Manuscripta Math., Vol. 115 (2004), 199--205. 


\bibitem{LS} Lytchak, A.; Stadler, S.,  \textit{Conformal deformations of CAT(0) spaces.} Math. Ann. 373 (2019), no. 1-2, 155--163. 


\bibitem{LW}Lytchak, A.; Wenger, S. ``Area minimizing discs in metric spaces,'' Arch. Ration. Mech. Anal., 223(3):1123--1182, 2017.

\bibitem{LW3}
Lytchak, A.; Wenger, S.,
\textit{Intrinsic structure of minimal discs in metric spaces.} 
Geom. Topol. 22 (2018), no. 1, 591--644. 

\bibitem{LW5} Lytchak, A.; Wenger, S.,  
\textit{Isoperimetric  characterization  of  upper  curvature bounds.}
Acta Math. 221 (2018), no. 1, 159--202. 


\bibitem{mese} Mese, C.
\textit{The curvature of minimal surfaces in singular spaces.}
Comm. Anal. Geom. 9 (2001), no. 1, 3--34. 

\bibitem{moore}
Moore, R. L.,
\textit{Concerning upper semi-continuous collections of continua.}
Trans. Amer. Math. Soc. 27 no. 4 (1925), 416--428.

\bibitem{petrunin-metric-min} Petrunin, A.,
\textit{Metric minimizing surfaces.}
Electron. Res. Announc. Amer. Math. Soc. 5 (1999), 47--54. 

\bibitem{petrunin-metric-min-correction} Petrunin, A.
\textit{Correction to: Metric minimizing surfaces.}
Electron. Res. Announc. Math. Sci. 25 (2018), 96. 

\bibitem{petrunin-intrinisic} Petrunin, A.,
\textit{Intrinsic isometries in Euclidean space.}
St. Petersburg Math. J. 22 (2011), no. 5, 803--812.

\bibitem{petrunin-stadler} Petrunin, A.; Stadler, S.,
\textit{Monotonicity of saddle maps.} 
Geom. Dedicata 198 (2019), 181--188.

\bibitem{shefel-2D} 
\begin{otherlanguage}{russian}
Шефель, С. З.,
\textit{О седловых поверхностях ограниченной спрямляемой кривой.}
Доклады АН СССР, 162 (1965) №2, 
294---296.
\end{otherlanguage}

\bibitem{shefel-3D} 
\begin{otherlanguage}{russian}
Шефель, С. З., 
\textit{О внутренней геометрии седловых поверхностей.}
Сибирский математический журнал, 5 (1964), 1382---1396.
\end{otherlanguage}

\bibitem{St} Stadler, S. \textit{The structure of minimal surfaces in CAT(0) spaces},
J. Eur. Math. Soc. (JEMS) 23 (2021), no. 11, 3521–3554.
  
\end{thebibliography}
\end{document}